%% file: template.tex
\documentclass[12pt]{article}
\usepackage{ulem}
\usepackage{amsmath}
\usepackage{mathrsfs}

\usepackage{graphicx} % 
\usepackage{caption}

\PassOptionsToPackage{hyphens}{url}

\usepackage[hyphens]{url}
\usepackage{hyperref}
\input{prefix}%
\newcommand{\bbS}{\mathbb{S}}
\newcommand{\barC}{\bar{\mathbb{C}}}

\begin{document}

\title{Two families of reducible spherical conical metrics}

\author{%
   Haoran Wu %
   \HaoranThanks{}%
   \and%
   Xuwen Zhu %
   \XuwenThanks{}%
}

\date{}

\maketitle

\begin{abstract}
We  analyze a 1-parameter family of heart shape and a 3-parameter family obtained by gluing three footballs, both of which are  examples of reducible spherical conical metrics. For these examples we verify the structure theorem given in~\cite{WeiWuXu2022} and show that such metrics naturally arise from Abelian differentials of the third kind. We then obtain the geometric decomposition using explicit metric and geodesic calculations. This offers new evidence for the interaction between the synthetic spherical geometry and the complex analytic structure of reducible conical metrics.
\end{abstract}

%%%%%%%%%%%%%%%%%%%%%%%%%%%%%%%%%%%%%%%%%%%%%%%%%%%%%%
%%%%%%%%%%%%%%%%%%%%%%%%%%%%%%%%%%%%%%%%%%%%%%%%%%%%%%

\section{Introduction}

%%% DELETE START 

%%%%%%%%%%%%%%%%%%%%%%%%%%%%%%%%%%%%%%%%%%%%%%%%
%%%
In this paper we study several examples of constant curvature conical metrics on the Riemann sphere. We call a metric $ds^{2}$ on a compact Riemann surface $\mathcal{M}$ a conical metric if it is smooth on $\mathcal{M}$ except at finitely many points $p_{1},\ldots,p_{N}$, where the conical singularity with angle $2\pi\beta_k$ at $p_k$ is locally given by $dr^2 + \beta_k^2 r^2  d\theta^2$. A conical metric with constant Gaussian curvature $K=-1,0$ or $1$ is called a hyperbolic, flat, or spherical conical metric, respectively. This paper focuses on the spherical case. 

There has been a long history and a lot of recent interest in investigating the existence and uniqueness of such metrics. 
Troyanov\cite{Troyanov1989} first demonstrated that a spherical conical metric with two singularities
of angles $\alpha$ and $\beta$ exists on the Riemann sphere $\bbS^2$ if and only if $\alpha=\beta$. In his later work~\cite{Troyanov1991}, he established the foundational existence results, providing necessary and sufficient conditions for spherical metrics with prescribed cone angles on $\mathbb{S}^2$ and more general compact Riemann surface $\mathcal{M}$. 
Luo and Tian\cite{LuoTian1992} later proved uniqueness for metrics on the sphere when all cone angles lie in $(0,2\pi)$. Chen and Li~\cite{ChenLi1991} obtained further necessary conditions under analytic constraints. Umehara and Yamada~\cite{UmeharaYamada2000} analyzed the three-cone-point case and introduced reducible spherical conical metrics (also called metrics with reducible monodromy), which was also studied by Eremenko~\cite{Eremenko2004}. In~\cite{BoLi} the angle condition for the existence of a reducible metric on an arbitrary Riemann surface was described. 
We refer to Chen-Wang-Wu-Xu~\cite{ChenWangWuXu2015} and Eremenko\cite{Eremenko2020} for detailed discussions of the reducible conical metrics.
One prime example of such reducible metrics is given by the "football", i.e. the spherical conical metric on $\mathbb{S}^2$ with two antipodal points with equal cone angles. 
In~\cite{Zhu2020} the second author analyzed the heart shape, which is a sphere with four conical points obtained by gluing together two footballs, showing by elementary trigonometry arguments that the resulting family exhibits local rigidity. This gave an evidence of the analytic obstruction considered in works of Mazzeo and the second author~\cite{MazzeoZhu2020,MazzeoZhu2022}. In~\cite{karpukhinZhu} the authors studied the spectral obstruction condition using harmonic maps, and constructed many examples of reducible metrics with such conditions, including one obtained by gluing three footballs. 

In~\cite{ChenWangWuXu2015} a link between reducible metrics and meromorphic differentials on Riemann surfaces was established. 
More recently Wei–Wu–Xu~\cite{WeiWuXu2022} developed a framework for describing reducible spherical conical metrics in terms of a character 1-form $\omega$, whose zeros and poles encodes cone angles and positions. Their structure theorem shows that every reducible spherical metric decomposes canonically into pieces isometric to cut-open footballs.

In this paper we will study several examples in detail, building on the work of Wei-Wu-Xu.
We  analyze a 1-parameter family of heart shape and a 3-parameter family obtained by gluing three footballs, which include examples studied in~\cite{Zhu2020} and~\cite{karpukhinZhu}. For these examples we verify the structure theorem given in~\cite{WeiWuXu2022} and show that such metrics naturally arise from Abelian differentials of the third kind, and then study their rigidity and geometric decomposition using explicit metric and geodesic calculations.

To state our result, we first recall a result on character 1-forms for reducible metrics from \cite{WeiWuXu2022}.
\begin{theorem}\label{thm1} \textbf{\upshape (\cite{WeiWuXu2022})} 
Suppose $ds^2$ is a reducible spherical conical metric on a compact Riemann surface $\mathcal{M}$, then there exists an Abelian differential $\omega$ of the third kind on $\mathcal{M}$ with
the following properties:
\begin{enumerate}
    \item All residues of $\omega$ are nonzero and real;
    \item $\omega$ has an exact real part outside its poles.
\end{enumerate}
If $p$ is a zero of $\omega$, then the cone angle of $ds^2$ at $p$ is $2\pi \cdot ( ord_p( \omega ) + 1)$. If $p$ is a pole for $\omega$, then the cone angle at $p$ is given by $2\pi\cdot|\mathrm{Res}_p(\omega)|$ unless $\mathrm{Res}_p(\omega)=\pm1$ (in this case $p$ is a smooth point of $ds^2$).  Furthermore, there exists a smooth
surjective function $\Phi:\mathcal{M}\setminus\{\textit{poles of }\omega\}\to(0,4)$  that can be continuously extended to  $\mathcal{M}$ and satisfies
\begin{equation}
    \frac{4\mathrm{d}\Phi}{\Phi(4-\Phi)}=\omega+\overline{\omega}. \label{differential equation}
\end{equation}
And the metric $ds^2$ is given by
\begin{equation}
    \mathrm{d}s^2=\frac{\Phi(4-\Phi)}{4}\omega\overline{\omega}. \label{ds2}
\end{equation}
Conversely, given an Abelian differential $\omega$ of the third kind on a compact Riemann surface $\mathcal{M}$ with properties 1 and 2, there exists a reducible spherical conical metric on $\mathcal{M}$ by~\eqref{ds2} with cone point information prescribed by $\omega$. 
\end{theorem}

%\begin{remark}
%    On a compact Riemann surface, an Abelian differential is referred to a meromorphic 1-form (see \cite{Spring1957}).
%\end{remark}
In this paper we will compute explicitly the 1-form $\omega$ in the two-football case and the three-football case. To streamline our discussion, we adopt the following notation:  we denote $\mathbb{S}_{\alpha_1,\ldots,\alpha_I}^2$ as the 2–sphere equipped with a spherical conical metric with cone angles $2\pi\alpha_1,\ldots,2\pi\alpha_I$. When only two angles appear and their distance is $\pi$, the surface is a football, written as $\bbS_{\alpha,\alpha}^{2}$.

We summarize our results as follows:
 \begin{theorem}\label{thm:heart}
      Let $ds^2 = \bbS^2_{2,\beta,\gamma}$ be the heart shape given by a spherical metric on $\mathbb{S}^2$ with angles $ 4\pi, 2\pi \beta, 2\pi \gamma$ with $\beta+\gamma=1$  (see Figure~\ref{fig:heartshape_decomposition}). Up to a change of coordinate $z \to pz, \ p\in\mathbb{C}\setminus\{0\}$, the corresponding Abelian differential of 3rd kind $\omega$ on $\mathbb{C}\cup\{\infty\}$ is given by
   \begin{equation}
    \omega=\frac{z}{(z-1)(z+\frac{\gamma}{\beta})}dz.
   \end{equation}  
   And there is a real one-parameter family of such metrics given by
   \begin{equation}
    ds^{2}=\frac{4e^{c}|z-1|^{-2(\beta-1)}|z+\frac{\gamma}{\beta}|^{-2(\gamma-1)}|z|^{2}}{[1+e^{c}|z-1|^{2(\beta-1)}|z+\frac{\gamma}{\beta}|^{2(\gamma-1)}]^{2}} |dz|^2. \label{ds2heart}
\end{equation}  
 \end{theorem}
\begin{remark}
    The poles of $\omega$ are given by $1, -\frac{\beta}{\gamma}, \infty$ with $Res_{1}(\omega)=\beta$, $Res_{-\frac{\gamma}{\beta}}(\omega)=\gamma$, and $Res_\infty(\omega)=-1$,  corresponding to the cone points of angles $2\pi\beta, 2\pi\gamma$ and a smooth point. And $\omega$ has one simple zero at 0 corresponding to a cone point with angle $4\pi$.
\end{remark} 

\begin{theorem}\label{thm:3football}
Let $ds^2$ be a spherical metric on $\mathbb{S}^2$ with angles $4\pi,4\pi,2\pi\beta, 2\pi(\alpha+\beta), 2\pi(\alpha+\gamma), 2\pi\gamma$ such that $\beta, \gamma, \alpha+\beta, \alpha+\gamma \notin \mathbb{Z}$. The corresponding Abelian differential of 3rd kind $\omega$ on $\mathbb{C}\cup\{\infty\}$ is given by
\begin{equation}
\omega=\left(\frac{-\beta}{z-P_{\beta}}+\frac{\alpha+\beta}{z-P_{\alpha}}+\frac{\gamma}{z-P_\gamma}\right)dz.
\end{equation}
where $P_\alpha, P_\beta, P_\gamma \in \barC\setminus\{0,1,\infty\}$ satisfy the following two equations:
\begin{equation} \label{g}
\begin{aligned} 
    &\beta P_\alpha P_\gamma+(\alpha+\beta)P_\beta P_\gamma+\gamma P_\beta P_\alpha=0 \\
    &-\beta(1-P_\alpha)(1-P_\gamma)+(\alpha+\beta)(1-P_\beta)(1-P_\gamma)+\gamma(1-P_\beta)(1-P_\alpha)=0.
\end{aligned}
\end{equation}
The set of $\{P_\alpha, P_\beta, P_\gamma\}\in (\barC\setminus\{0,1,\infty\})^3$ forms a (complex) 1-dimensional algebraic subvariety. There is a real 3-parameter family of metrics given by
\begin{equation}
ds^2=\frac{4c^2\beta^2(\alpha+\beta)^2\gamma^2|(z-P_{\beta})^{-\beta-1}(z-P_{\alpha})^{\alpha+\beta-1}(z-P_{\gamma})^{\gamma-1}|^2}{\left(1+|c(z-P_{\beta})^{-\beta}(z-P_{\alpha})^{\alpha+\beta}(z-P_{\gamma})^{\gamma}|^2\right)^2}dz^2.
\end{equation}
 \end{theorem}

\begin{remark}\remlab{remark1.1}
    In the three-football case with the angle combination given above, there is a real 3-parameter family of metrics, see Figure~\ref{fig:threefootball} and~\ref{fig:threefootball_v2} for two distinct configurations and Figure~\ref{fig:threefootball_v2} for a generic configuration with marked parameters. Roughly speaking, the parameters correspond to the choice of a point in the 1-dimensional algebraic subvariety mentioned above and a conformal scaling factor in the metric. We will show later that one way to see the three parameters is through the lengths of  geodesics. 
    %The way to distinguish these two cases is by computing the lengths of geodesics between cone points. From \cite[Example 7.4.]{WeiWuXu2022}, we know that in case 2 (see Figure~\ref{fig:threefootball_v2}) the two saddle points (0 and 1) are located on a common geodesic that directly links a minimum point to a maximum point of $\Phi$ in \eqref{differential equation} on the sphere, and such geodesic has the length of $\pi$.
\end{remark}

 The work of Tahar~\cite{Tahar2022} described how reducible spherical metrics can be decomposed by cutting along finitely many geodesics joining cone points. In~\cite{WeiWuXu2022} the above mentioned Theorem~\ref{thm1} was also used to prove a structure theorem about decomposition of reducible metrics into footballs. In the examples we study in this paper, we will use explicit calculations of geodesics to verify that the heart-shaped metric $\bbS^2_{2,\beta,\gamma}$ admits a decomposition into two footballs $\bbS^2_{\beta, \beta}$ and $\bbS^2_{\gamma,\gamma}$ sharing a common boundary (see Figure~\ref{fig:heartshape_decomposition}), and the  metric $\bbS^2_{2, 2, \beta, \alpha+\beta, \alpha+\gamma, \gamma}$ can be cut into three footballs $\bbS^2_{\alpha,\alpha}, \bbS^2_{\beta,\beta}$ and $\bbS^2_{\gamma,\gamma}$ (see Figure~\ref{fig:threefootball} and~\ref{fig:threefootball_v2} for two examples). 

We also mention that recently Chen et al.\cite{ChenEtAl2025} developed a detailed analysis of isoresidual strata for meromorphic and k-differentials and presented an application to cone spherical metrics with dihedral monodromy which is a generalization of reducible monodromy. See also~\cite{GendronTahar2023}. One expects similar calculations can be applied to those cases as well.

This paper is divided into two sections. Section 2 describes the calculation on the heart shape, where we compute the meromorphic differential and the spherical metric, then use the metric to determine geodesics between cone points and the decomposition. Section 3 works on the three-football case, where a similar calculation follows, and we will show the existence of a 3-parameter family of configurations using geodesic calculation.

\noindent \textbf{Acknowledgement:} We would like to thank Bin Xu for helpful comments. The second author is supported by NSF DMS-2305363.

%------------------------------------------------------------------
%------------------------------------------------------------------

\section{Heart shape}

\begin{figure}[htbp]
  \centering
  \includegraphics[width=1\textwidth]{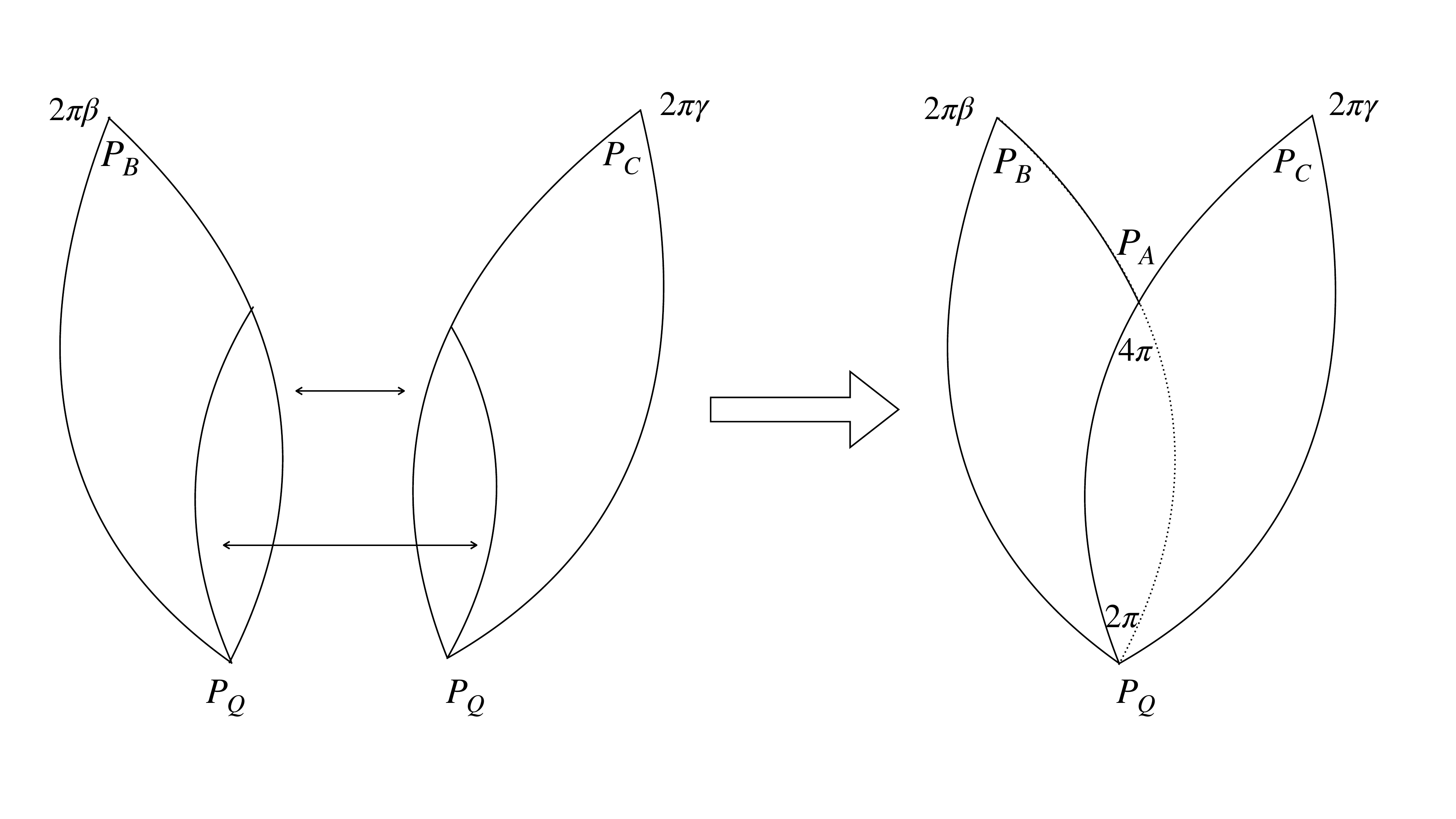} % 
  \caption{A picture of $\bbS_{\{2,\beta,\gamma\}}^2$ with $\beta+\gamma=1$, obtained from gluing two footballs $\bbS^2_{\beta,\beta}$ and $\bbS^2_{\gamma, \gamma}$}
  \label{fig:heartshape_decomposition}
\end{figure}

We will compute the differential $\omega$ of the heart shape, see  Figure \ref{fig:heartshape_decomposition}. This example was studied in~\cite[Figure 1] {Zhu2020} and \cite[Example 7.1]{WeiWuXu2022}, and we will use it to demonstrate the calculation and provide a more detailed description of the geometry. Let $\beta, \gamma \in (0,1)$ be two numbers satisfying $\beta+\gamma=1$, we consider a reducible spherical conical metric on $\mathbb{S}^{2}$ with three cone points $4\pi, 2\pi\beta, 2\pi\gamma$, denoted by $S_{\{2,\beta,\gamma\}}^2$. One way to construct such a metric is to glue two footballs, $S_{\{\beta,\beta\}}^2$ and $S_{\{\gamma,\gamma\}}^2$ along a cut on the meridian. The calculation below will show that this is only way to obtain such a metric.

\begin{remark}
    The above method can be used to construct a reducible spherical conical metric on any value of $\beta$ and $\gamma$,  see \cite[Example 7.2]{WeiWuXu2022}. We chose the condition $\beta+\gamma=1$ for the purpose of simplifying the calculation for $\omega$, but a similar calculation holds for arbitrary $\beta$ and $\gamma$.
\end{remark}

To find the meromorphic differential in Theorem~\ref{thm1},  we consider $\omega$ defined on $\bar{\mathbb{C}} = \mathbb{C}\cup \infty$ such that the conical points of angles $4\pi, 2\pi\beta, 2\pi\gamma$ are denoted by $P_A, P_B$ and $P_C \in \bar{\mathbb{C}}$. Recall if $p$ is a zero of $\omega$ then the cone angle is $2\pi(ord_p(\omega)+1)$; if $p$ is a pole of $\omega$, then the cone  angle is $\pm2\pi\cdot Res_p(w)$. Since $\beta$ and $\gamma$ are not integers, $P_B$, $P_C$ must be poles for $\omega$. On the other hand, the remaining cone point $P_A$ must be a zero, and we set  $P_A=0$ without loss of generality (otherwise one can perform a conformal transformation on $\bar{\mathbb{C}}$). We also denote $P_Q = \infty$. Since we have $ord_{P_A} = 1$, $Res_{P_C}=\gamma$, $Res_{P_B}=\beta$,  by Global Residue Theorem we know $P_Q$ is also a pole of $\omega$ and $Res_{P_Q}=-1$ which corresponds to a smooth point. Now we are ready to state the theorem:

 \begin{theorem}\thmlab{heartshape}
    (same as Theorem~\ref{thm:heart})  Let $ds^2$ be the heart shape given by a spherical metric on $\mathbb{S}^2$ with angles $2\pi \beta, 2\pi \gamma, 4\pi$ with $\beta+\gamma=1$  (see Figure~\ref{fig:heartshape_decomposition}). Up to a change of coordinate $z \to pz, \ p\in\mathbb{C}\setminus\{0\}$, the corresponding Abelian differential of 3rd kind $\omega$ on $\mathbb{C}\cup\{\infty\}$ is given by
   \begin{equation}
    \omega=\frac{z}{(z-1)(z+\frac{\gamma}{\beta})}dz
   \end{equation}
 And there is a real one-parameter family of such metrics given by
   \begin{equation}
    ds^{2}=\frac{4e^{c}|z-1|^{-2(\beta-1)}|z+\frac{\gamma}{\beta}|^{-2(\gamma-1)}|z|^{2}}{[1+e^{c}|z-1|^{2(\beta-1)}|z+\frac{\gamma}{\beta}|^{2(\gamma-1)}]^{2}} |dz|^2. \label{ds2heart}
\end{equation}   
 \end{theorem}

\noindent{\textit{Proof.}}
We will follow the proof of~\cite[lemma 3.1]{WeiWuXu2022}. Recall the general form of $\omega$ using the information of two poles at $P_B$ and $P_C$: 
\begin{equation}
    \omega=f(z) dz =(\frac{\beta}{z-P_{B}}+\frac{\gamma}{z-P_{C}})dz \label{omega}.
\end{equation}
In other words
$$\omega =\frac{\beta z-\beta P_{C}+z\gamma-\gamma P_{B}}{(z-P_{B})(z-P_{C})}dz.$$
There is a zero of order 1 at $z=0$, therefore the locations of the two poles need to satisfy
\begin{equation}
    P_{C}=-\frac{\gamma}{\beta} P_{B}. \label{bcrelation}
\end{equation}
Now we check the residue at $\infty$, which is given by $-\lim\limits_{|z|\to\infty}(\frac{\beta z}{z-P_B}+\frac{\gamma z}{z-P_C})=-(\beta+\gamma)=-1$, hence it matches the condition that $\infty$ is a smooth point. After we apply the relation (\ref{bcrelation}) to (\ref{omega}), we have
 $$w=(\frac{\beta}{z-P_{B}}+\frac{\gamma}{z+\frac{\gamma}{\beta}P_{B}})dz.$$
 Now we consider the change of coordinate $z \to pz, \ p\in \mathbb{C}$,
 $$\begin{aligned}\omega&=(\frac{\beta}{pz-P_{B}}+\frac{\gamma}{pz+\frac{\gamma}{\beta}P_{B}})pdz\\&=(\frac{\beta}{z-\frac{P_{B}}{p}}+\frac{\gamma}{z+\frac{\gamma}{\beta}\frac{P_{B}}{p}})dz\end{aligned}$$
 Set $p=P_B$, we have the desired form of $\omega$:
 $$\begin{aligned}\omega&=(\frac{\beta}{z-1}+\frac{\gamma}{z+\frac{\gamma}{\beta}})dz\\&=\frac{z}{(z-1)(z+\frac{\gamma}{\beta})}dz\end{aligned}.$$

With this differential form, we will now use \eqref{differential equation} to get $\Phi$ and use \eqref{ds2} to get the metric $ds^2$. Recall that since $\omega+\bar{\omega}$ is exact on $\mathcal{M}'=\mathcal{M}\setminus\{0,1,-\frac{\gamma}{\beta}, \infty\}$, there exists a smooth function $f$ defined on $\mathcal{M}'$ such that
$$ \frac{4\mathrm{d}\Phi}{\Phi(4-\Phi)}=\omega+\overline{\omega}=\mathrm{d}f.$$
Using
$$\frac{4\mathrm{d}\Phi}{\Phi(4-\Phi)}=\mathrm{d}\ln\frac\Phi{4-\Phi},$$
we have
$$\ln\frac\Phi{4-\Phi}=f+c$$
which implies
\begin{equation}
    \Phi=\frac{4e^{f+c}}{1+e^{f+c}}\in(0,4). \label{phi}
\end{equation}
In the heart shape case we get  
$$f(z)=\beta\ln|z-1|^{2}+\gamma\ln|z+\frac{\gamma}{\beta}|^{2}+c,$$ where $c$ is a real constant. Plugging into \eqref{phi}, we have 
$$\Phi=\frac{4e^{c}|z-1|^{2\beta}\cdot|z+\frac{\gamma}{\beta}|^{2\gamma}}{1+e^{c}|z-1|^{2\beta}\cdot|z+\frac{\gamma}{\beta}|^{2\gamma}}.$$
By \eqref{ds2}, we get the metric
\begin{equation}
    ds^{2}=\frac{4e^{c}|z-1|^{-2(\beta-1)}|z+\frac{\gamma}{\beta}|^{-2(\gamma-1)}|z|^{2}}{[1+e^{c}|z-1|^{2(\beta-1)}|z+\frac{\gamma}{\beta}|^{2(\gamma-1)}]^{2}} |dz|^2.
\end{equation}

\qed

\begin{figure}[htbp]
  \centering
  \includegraphics[width=0.8\textwidth]{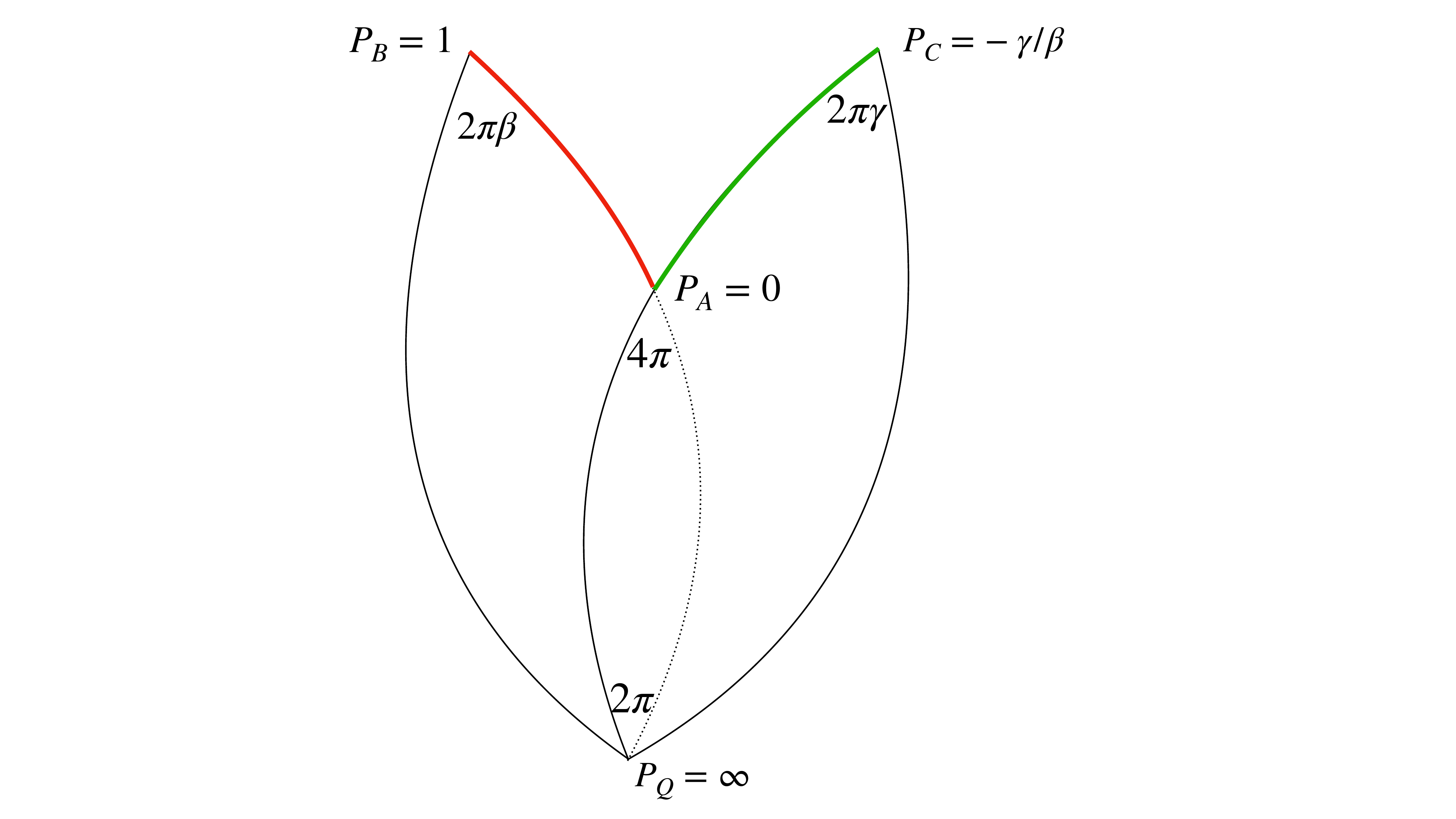} % 
  \caption{}
  \label{fig:heartshape_length}
\end{figure}

Now let us consider the geodesic from $0$ to $1$ on $\bar{\mathbb{C}}$, corresponding to the red line in Figure \ref{fig:heartshape_length}, and from $0$ to $-\gamma/\beta$, which is the green line in Figure \ref{fig:heartshape_length}. 
These two geodesics connect $P_A$ to $P_B$ and $P_C$. 
We claim that their lengths are the same and depend only on the constant $c$ in~\eqref{ds2heart}.  Let $F(z)=e^{\frac{c}{2}}(z-1)^{\beta}(z+\frac{\gamma}{\beta})^{\gamma}$ be a multi-valued meromorphic function on $\bar{\mathbb{C}}$ we have 
$$ds^{2}=\frac{4|F^{\prime}(z)|^{2}}{(1+|F(z)|^{2})^{2}}|dz|^{2}.$$
Recall the standard round metric on the sphere
$$ds^{2}_{\mathbb{S}^2}=\frac{4|dw|^{2}}{(1+|w|^{2})^{2}}, \ w\in \bar{\mathbb{C}},$$
then the heart shape metric can be obtained as the pullback of the standard sphere metric by $F$ (which is the developing map \cite{ChenWangWuXu2015}), such that
$$ds^2=F^*ds^{2}_{\mathbb{S}^2}=\frac{4|F^{\prime}(z)|^{2}}{(1+|F(z)|^{2})^{2}}|dz|^{2}.$$
The geodesics of the pulled back metrics are (locally) preimages under $F$ of great circles on the sphere. Consider the images of $P_A$ and $P_B$ under $F$ ($F$ is multivalued but we can pick any of the images for $F(0)$):
$$w_0:=F(0)=e^{\frac{c}{2}}(-1)^{\beta}(\frac{\gamma}{\beta})^{\gamma},\ w_1:=F(1)=0.$$
The geodesic between $w_0$ and $w_1$ is the great circle between the south pole and $\omega_0$. Using stereographic projection, we see that the geodesic is a straight line through the origin
$$L=\{tw_0:t\in [0,1]\}.$$
The geodesic in $ds^2$ is its preimage: $\gamma=F^{-1}(L)=\{z:F(z)\in L\}$. That is, for all $z$ we should have 
$$\frac{F(z)}{F(0)}=\frac{(z-1)^\beta\left(z+\frac{\gamma}{\beta}\right)^\gamma}{(-1)^\beta\left(\frac{\gamma}{\beta}\right)^\gamma}\in [0,1].$$
Note that locally away from  all the cone points $F$ is a diffeomorphism, so we can choose the branch of $F(z)$ based on the choice of $F(0)$. Hence the geodesic is
$$\gamma=\left\{z\in\mathbb{C}:\frac{(z-1)^\beta\left(z+\frac{\gamma}{\beta}\right)^\gamma}{(-1)^\beta\left(\frac{\gamma}{\beta}\right)^\gamma}\in [0,1]\right\}.$$
Now let us consider the length of this geodesic segment 
$$L(\gamma)=\int_{\gamma}ds=2\int_{\gamma}\frac{|F^{\prime}(z)|}{1+|F(z)|^{2}}|dz|$$
Using a change of variables,  we set $t=F(z)/F(0)\in [0,1]$, $z=F^{-1}(tF(0))$. Note that $dz=\frac{F(0)dt}{F^{\prime}(z)}$, we obtain
$$L(0,1)=2\left|\int_{t=0}^{t=1}\frac{|F^{\prime}(z)|}{1+|F(0)|^2t^2}\cdot\frac{|F(0)|dt}{|F^{\prime}(z)|}\right|=2\left|\int_0^1\frac{|F(0)|dt}{1+t^2|F(0)|^2}\right| = 2\left|\int_0^{|F(0)|} \frac{ds}{1+s^2}\right|=2\arctan(|w_0|) .$$
Similar arguments apply to lengths from $0$ to $-\gamma/\beta$, and from $0$ to $\infty$, and we have
$$L(0,-\gamma/\beta)=2\left|\int_{|F(-\gamma/\beta)|}^{|F(0)|}\frac{ds}{1+s^2}\right|=2\arctan(|w_0|), $$
$$L(0,\infty)=2\left|\int_{|F(\infty)|}^{|F(0)|}\frac{ds}{1+s^2}\right|=\pi-2 \arctan(|w_0|). $$
Note there are two geodesics from $0$ to $\infty$ because of the choice of branches, but the length $L(0,\infty)$ is the same.

\begin{figure}[htbp]
  \centering
  \includegraphics[width=0.8\textwidth]{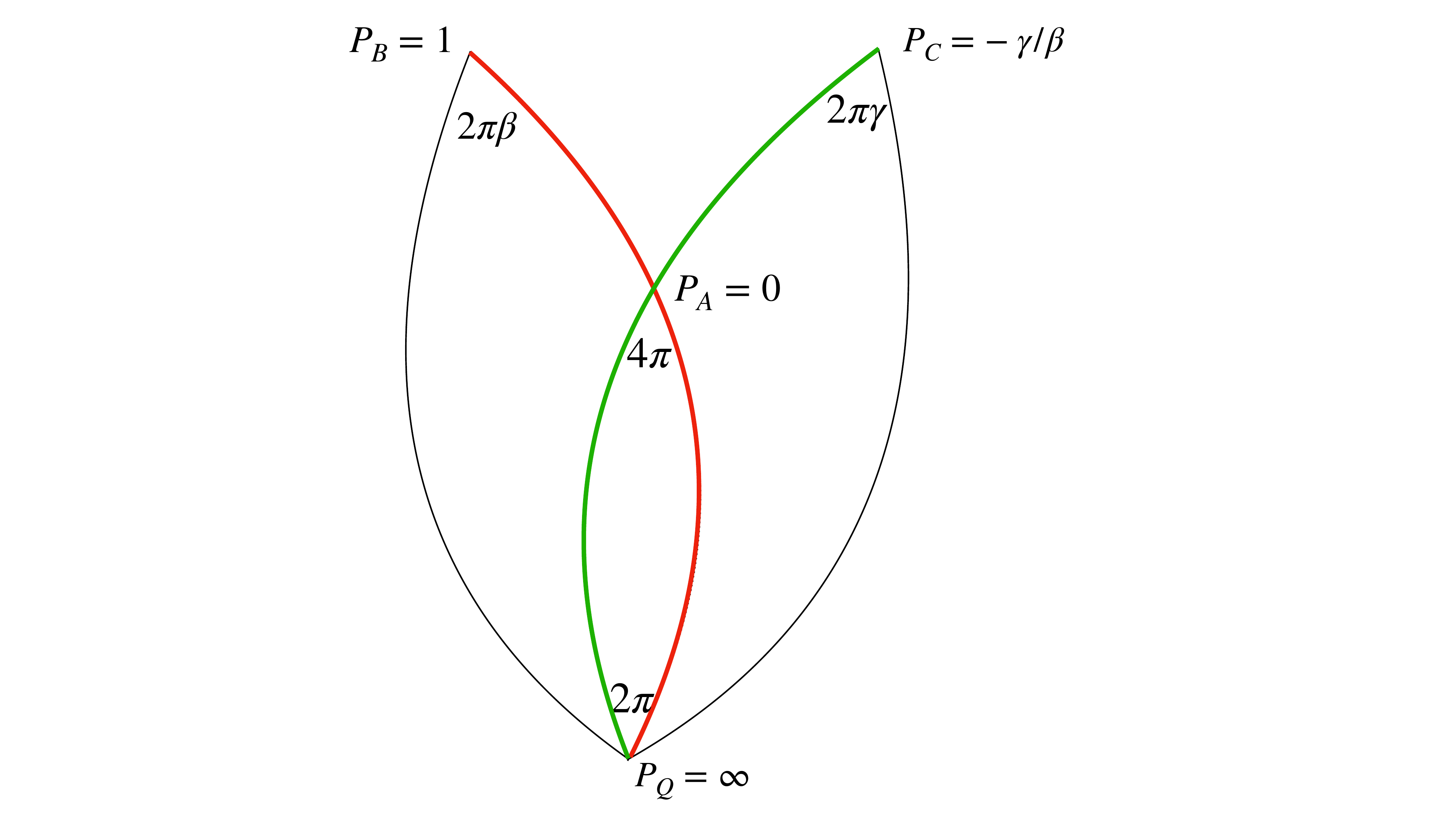} % 
  \caption{Geodesics from $1$ to $\infty$ and from $-\frac{\gamma}{\beta}$ to $\infty$. They both go through 0 and our calculation shows that the lengths are both $\pi$.}
  \label{fig:heartshape_totallength}
\end{figure}

Note that we have $L(0,1)=L(0,-\frac{\gamma}{\beta})$ which only depends on the real parameter $c$. And we also have $L(0,1) + L(0, \infty) = \pi$ and $L(0,-\gamma/\beta) + L(0, \infty) = \pi$, representing the red line and green geodesics in Figure \ref{fig:heartshape_totallength} respectively. In this way we recover the gluing procedure of the heart shape.

\section{Three-football shape}

\begin{figure}[htbp]
  \centering  \includegraphics[width=0.6\textwidth]{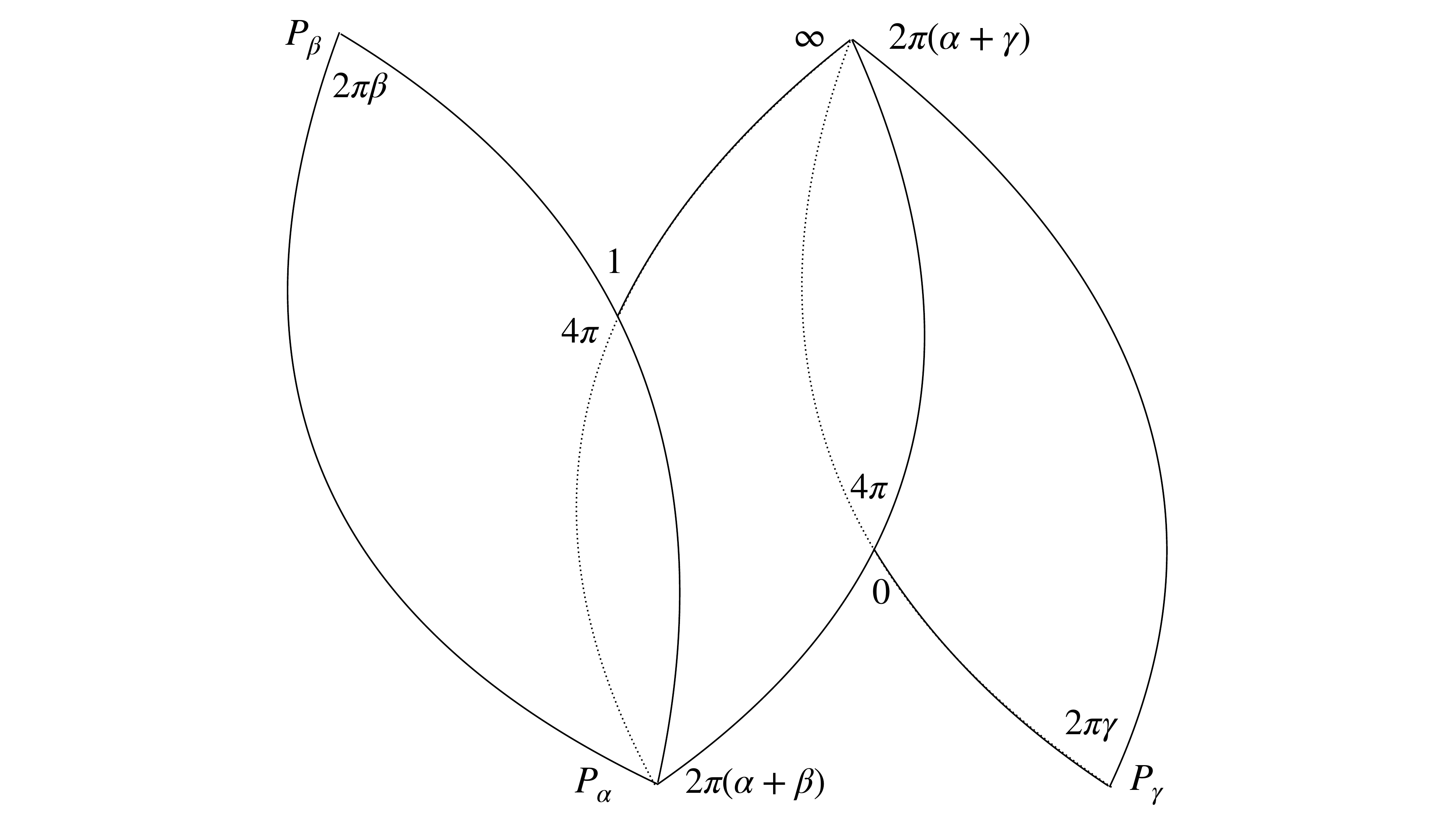} % 
  \caption{One example of $\bbS^2_{2,2,\beta,\alpha+\beta, \alpha+\gamma,\gamma}$ obtained by gluing three footballs $\bbS^2_{\beta, \beta}, \bbS^2_{\alpha, \alpha}, \bbS^2_{\gamma, \gamma}$}
  \label{fig:threefootball}
\end{figure}

Now we will consider a three-football case (see Figure \ref{fig:threefootball} for one of the examples). One can also see~\cite[Example 7.4. (4)]{WeiWuXu2022}. 
We will be looking for metrics $\mathbb{S}^2_{2,2,\beta, \alpha+\beta,\alpha+\gamma, \gamma}$. We will show that there is a real 3-parameter family (see Figure \ref{fig:family of threefootbal}) and one can recover the parameters by the lengths of geodesics. 

To set up the calculation for the 1-form $\omega$, we put two simple zeros at 0 and 1, corresponding to the conical points with angles $4\pi$. Let $P_\gamma$ be a pole with $Res_{P_\gamma} = \gamma$, $P_\beta$ be a pole with $Res_{P_\beta} = \beta$ and $P_\alpha$ be a pole with $Res_{P_\alpha} = \alpha+\beta$, where the corresponding angles are $2\pi\gamma, 2\pi\beta$ and $2\pi(\alpha+\beta)$. By the Global Residue Theorem, we can get $Res_\infty=-\alpha-\gamma$, and the cone angle there is $2\pi(\alpha+\gamma)$. Note that, with this given angle set with all the non-integer assumptions, the residue distribution we describe here is the only possible combination up to a uniform factor of $\pm1$. Now we can state the theorem:
\begin{theorem}\thmlab{threefootball}
(same as Theorem~\ref{thm:3football}) Let $ds^2$ be a spherical metric on $\mathbb{S}^2$ with angles $4\pi,4\pi,2\pi\beta, 2\pi(\alpha+\beta), 2\pi(\alpha+\gamma), 2\pi\gamma$ such that $\beta, \gamma, \alpha+\beta, \alpha+\gamma \notin \mathbb{Z}$. The corresponding Abelian differential of 3rd kind $\omega$ on $\mathbb{C}\cup\{\infty\}$ is given by
\begin{equation}
\omega=\left(\frac{-\beta}{z-P_{\beta}}+\frac{\alpha+\beta}{z-P_{\alpha}}+\frac{\gamma}{z-P_\gamma}\right)dz.
\end{equation}
where $P_\alpha, P_\beta, P_\gamma \in \barC\setminus\{0,1,\infty\}$ satisfy the following two equations:
\begin{equation} \label{g}
\begin{aligned} 
    &\beta P_\alpha P_\gamma+(\alpha+\beta)P_\beta P_\gamma+\gamma P_\beta P_\alpha=0 \\
    &-\beta(1-P_\alpha)(1-P_\gamma)+(\alpha+\beta)(1-P_\beta)(1-P_\gamma)+\gamma(1-P_\beta)(1-P_\alpha)=0.
\end{aligned}
\end{equation}
and the set of $\{P_\alpha, P_\beta, P_\gamma\}\in (\barC\setminus\{0,1,\infty\})^3$ forms a (complex) 1-dimensional algebraic subvariety. There is a real 3-parameter family of metrics given by
\begin{equation}\label{eq:ds3football}
ds^2=\frac{4c^2\beta^2(\alpha+\beta)^2\gamma^2|(z-P_{\beta})^{-\beta-1}(z-P_{\alpha})^{\alpha+\beta-1}(z-P_{\gamma})^{\gamma-1}|^2}{\left(1+|c(z-P_{\beta})^{-\beta}(z-P_{\alpha})^{\alpha+\beta}(z-P_{\gamma})^{\gamma}|^2\right)^2}dz^2.
\end{equation}
 \end{theorem}

\noindent\textit{Proof.} From the assumption of poles, we have the following
\begin{equation}\label{eq:omega}
    \omega=(\frac{-\beta}{z-P_{\beta}}+\frac{\alpha+\beta}{z-P_{\alpha}}+\frac{\gamma}{z-P_{\gamma}})dz=f(z)dz
\end{equation}
Denote $f(z) =\frac{g(z)}{h(z)}$ after taking the common denominator, we have  
$$g(z)=-\beta(z-P_\alpha)(z-P_\gamma)+(\alpha+\beta)(z-P_\beta)(z-P_\gamma)+\gamma(z-P_\beta)(z-P_\alpha)$$
Since $g(1)=g(0)=0$ we get the condition~\eqref{g}. $\{P_\alpha, P_\beta, P_\gamma\}$ are solutions to these two quadratic equations, therefore form a (complex) 1-dimensional algebraic variety. 

One could also start with the following form of $\omega$ taking into account the simple zeros and the positions of poles:
\begin{equation}
\omega=C\frac{z(z-1)}{(z-P_\beta)(z-P_\gamma)(z-P_\alpha)}dz, \ C\in \mathbb{C}\setminus \{0\}.
\end{equation}
Using the residue information at $P_\alpha, P_\beta, P_\gamma$ we get three quadratic equations for $\{C, P_\alpha, P_\beta, P_\gamma\} \in  \mathbb{C}\setminus \{0\} \times (\barC\setminus\{0,1,\infty\})^3$ which cuts out a 1-dimensional subvariety. 

Now we need to consider $\Phi$, by the same process as we did in the previous section we have:
\begin{equation}
    \Phi=\frac{4e^{f+A_{0}}}{1+e^{f+A_{0}}}
\end{equation}
where $f=-\beta\ln|z-P_\beta|^2+(\alpha+\beta)\ln|z-P_\alpha|^2+\gamma\ln|z-P_\gamma|^2$. 
%So we have
%\begin{equation}
%    \begin{aligned}
%   \frac{d\Phi}{df}&=\frac{4e^{f+A_{0}}(1+e^{f+A_{0}})-4e^{f+A_{0}}\cdot e^{f+A_{0}}}{(1+e^{f+A_{0}})^{2}}\\&
%    =\frac{4e^{f+A_{0}}}{(1+e^{f+A_{0}})^{2}}>0
%\end{aligned}     
%\end{equation}
%So $\Phi$ is a monotonically increasing function with respect to $f$. When $f\rightarrow-\infty$, $\Phi$ reaches its minimum, when $f\rightarrow\infty$ $\Phi$ reaches its maximum. When we consider the geodesic that travel from $P_\alpha$ to $\infty$ passing through the saddle points $0$ and $1$, it is easy to get that $P_\alpha$ is the minimum and $\infty$ is the maximum. So we have 
%\begin{equation}
%L_{total}=L(P_{\alpha},1)+L(1,0)+L(0,P_{\infty})
%\end{equation}
%By above argument and \remref{remark1.1}, we have that in Case 1, we have such length larger than $\pi$; and in  Case 2, we have such length equal to $\pi$, see the red line in Figure \ref{fig:threefootball_totallength}.

%\begin{figure}[htbp]
  %\centering  \includegraphics[width=0.7\textwidth]%{threefootball_totallength.pdf} % 
  %\caption{The two extreme cases can be distinguished by the %geodesic going from $P_\alpha$ to $\infty$.}
%\label{fig:threefootball_totallength}
%\end{figure}
We compute the developing map as
\begin{equation}\label{eq:F}
F=c(z-P_{\beta})^{-\beta}(z-P_{\alpha})^{\alpha+\beta}(z-P_{\gamma})^{\gamma}, c>0
\end{equation}
then by the same argument we did in heart shape case we have
\begin{equation}\label{eq:ds}
ds^2 = \frac{4|F^{\prime}(z)|^{2}}{(1+|F(z)|^{2})^{2}}|dz|^{2}
\end{equation}
Plugging $F(z)$ into~\eqref{eq:ds} we get the metric form in~\eqref{eq:ds3football}.
\qed

It remains to describe the geometry of the corresponding conical metric. To have an explicit calculation, we will focus on the solutions for a special case.
Recall the equations for pole positions
\begin{align}
    &g(0)-\beta P_\alpha P_\gamma+(\alpha+\beta)P_\beta P_\gamma+\gamma P_\beta P_\alpha=0 \label{g(0)}\\&g(1)=-\beta(1-P_\alpha)(1-P_\gamma)+(\alpha+\beta)(1-P_\beta)(1-P_\gamma)+\gamma(1-P_\beta)(1-P_\alpha)=0 \label{g(1)}
\end{align}
By (\ref{g(0)}), we obtain the following relation
\begin{equation}
P_{\alpha}=\frac{(\alpha+\beta)P_{\beta}P_{\gamma}}{\beta P_{\gamma}-\gamma P_{\beta}} \label{Palpha}
\end{equation}
If we plug (\ref{Palpha}) into (\ref{g(1)}), we have a quadratic form relating $P_\beta$ and $P_\gamma$
$$-\alpha\beta P_{\gamma}^{2}+(\alpha\beta+\gamma\beta-2\gamma\beta P_{\beta})P_{\gamma}+(-\alpha\gamma P_{\beta}+\alpha\gamma P_{\beta}^{2}+\gamma\beta P_{\beta}^{2}-\gamma^{2}P_{\beta}+\gamma^{2}P_{\beta}^{2})=0$$
Computing the determinant of the quadratic form $a P_\gamma^2 + b P_\gamma +c$, we have 
\begin{equation}
    \begin{aligned}b^{2}-4ac=&(-4\alpha\gamma^{2}\beta+4\alpha\gamma\beta^{2}+4\gamma^{2}\beta^{2}+4\alpha^{2}\gamma\beta)P_{\beta}^{2}+\\&(-4\alpha\gamma\beta^{2}-4\gamma^{2}\beta^{2}-4\alpha\gamma^{2}\beta)P_{\beta}+\\&\alpha^{2}\beta^{2}+2\alpha\gamma\beta^{2}+\gamma^{2}\beta^{2}\end{aligned} \label{b2-4ac}
\end{equation}
We will assume $\alpha=\gamma,\beta=\frac{(1+\sqrt{2})}{2}\alpha$, then (\ref{b2-4ac}) becomes  $b^2-4ac=(2\sqrt{2}\alpha \beta P_\beta-2\alpha\beta)^{2}$ hence we can simplify the quadratic form above. Assuming $P_\beta\neq \frac{1}{\sqrt{2}}$ we have 
$$P_\gamma=1-P_{\beta}\pm(-\sqrt{2}P_{\beta}+1).$$
Choose one of the solutions
\begin{equation}
    P_{\gamma}^{-}=(\sqrt{2}-1)P_{\beta} \label{tworoots}
    %, \ P_{\gamma}^{+}=2-(1+\sqrt{2})P_{\beta}. 
\end{equation} 
Plugging into~\eqref{Palpha} we can get
$$
P_\alpha^-=(1-2\sqrt{2}) P_\beta.
$$
Plugging into the expression of $\omega$ in~\eqref{eq:omega} we get 
\begin{equation}\label{eq:omega2}
\omega=\left(\frac{-\beta}{z-P_{\beta}}+\frac{\alpha+\beta}{z-(1-2\sqrt{2}) P_\beta}+\frac{\gamma}{z-(\sqrt{2}-1)P_\beta}\right)dz
\end{equation}
From this form one can see that locally the  space of $\omega$ is parametrized by $P_\beta \in U \subset \mathbb{C}$ where $U$ is a small open set outside $0,1,\infty$ or $\frac{1}{\sqrt{2}}$.

Putting all that into $F(z)$ in~\eqref{eq:F} we get
$$F(z)=c(z-P_\beta)^{-\beta} (z-(1-2\sqrt{2})P_\beta)^{\alpha+\beta} (z-(\sqrt{2}-1)P_\beta)^\gamma.$$
The geodesic length between $a,b \in \barC$ is
\begin{equation}
L(a,b)=2\left|\int_{|F(a)|}^{|F(b)|}\frac{1}{1+t^{2}}dt\right|=2|\arctan |F(b)|-\arctan |F(a)|| \label{length}
\end{equation}
if $a$ or $b$ is  $P_\alpha$ or $\infty$. (By the same argument as before, these geodesics are preimages of $F$ of geodesics going through 0 or $\infty$ hence can be described as straight lines in $\mathbb{C}$.)
Therefore we have
\begin{equation}\label{eq:geodesic}
    L(1, P_\infty)=\pi-2\arctan |F(1)|, \ L(0,P_{\infty}) = \pi-2\arctan |F(0|.
\end{equation}
Note that although $F(z)$ is multivalued, the geodesic lengths in~\eqref{eq:geodesic} are well-defined. We can see that the two parameters $c\in \mathbb{R}\setminus 0$ and $P_\beta \in U\subset \mathbb{C}$ determine the lengths $\ell_1 = L(1, P_\infty)$ and $\ell_2=L(0,P_{\infty})$, representing the red 
and green geodesics in Figure \ref{fig:family of threefootbal} respectively.

%When the saddle points are not located on a common geodesic that directly links a minimum point to a maximum point of $\Phi$,
%   \begin{equation}\label{eq:omega1}
 %  \omega=\left(\frac{-\beta}{z-P_{\beta}}+\frac{\alpha+\beta}{z-P_{\alpha}}+\frac{\gamma}{z-2+%(1+\sqrt{2})P_\beta}\right)dz
%   \end{equation}

\begin{figure}[htbp]
  \centering  \includegraphics[width=0.7\textwidth]{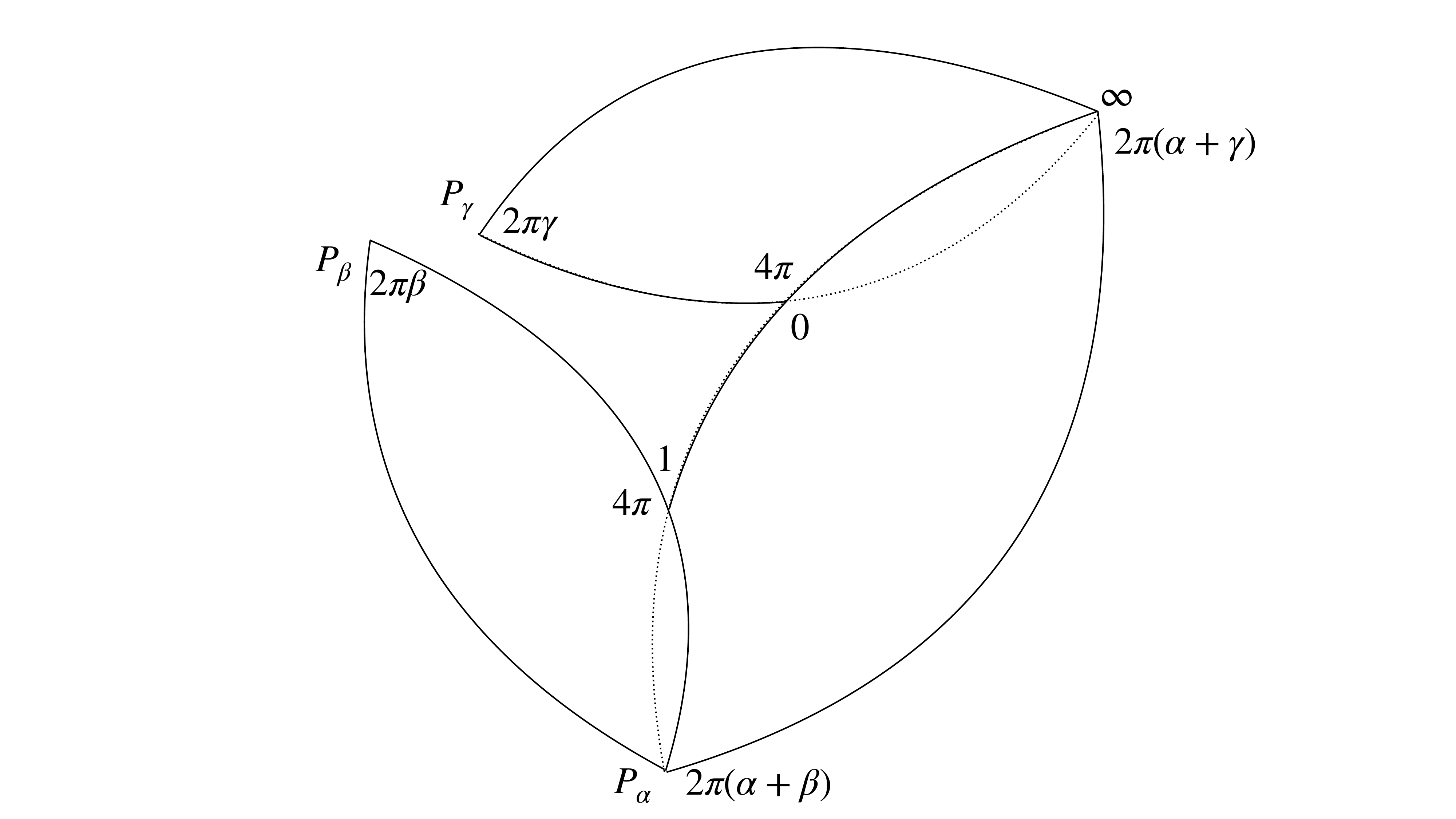} % 
  \caption{Another example of $\bbS^2_{2,2,\beta,\alpha+\beta,\alpha+\gamma, \gamma}$}
  \label{fig:threefootball_v2}
\end{figure}

\begin{figure}[htbp]
  \centering  \includegraphics[width=0.9\textwidth]{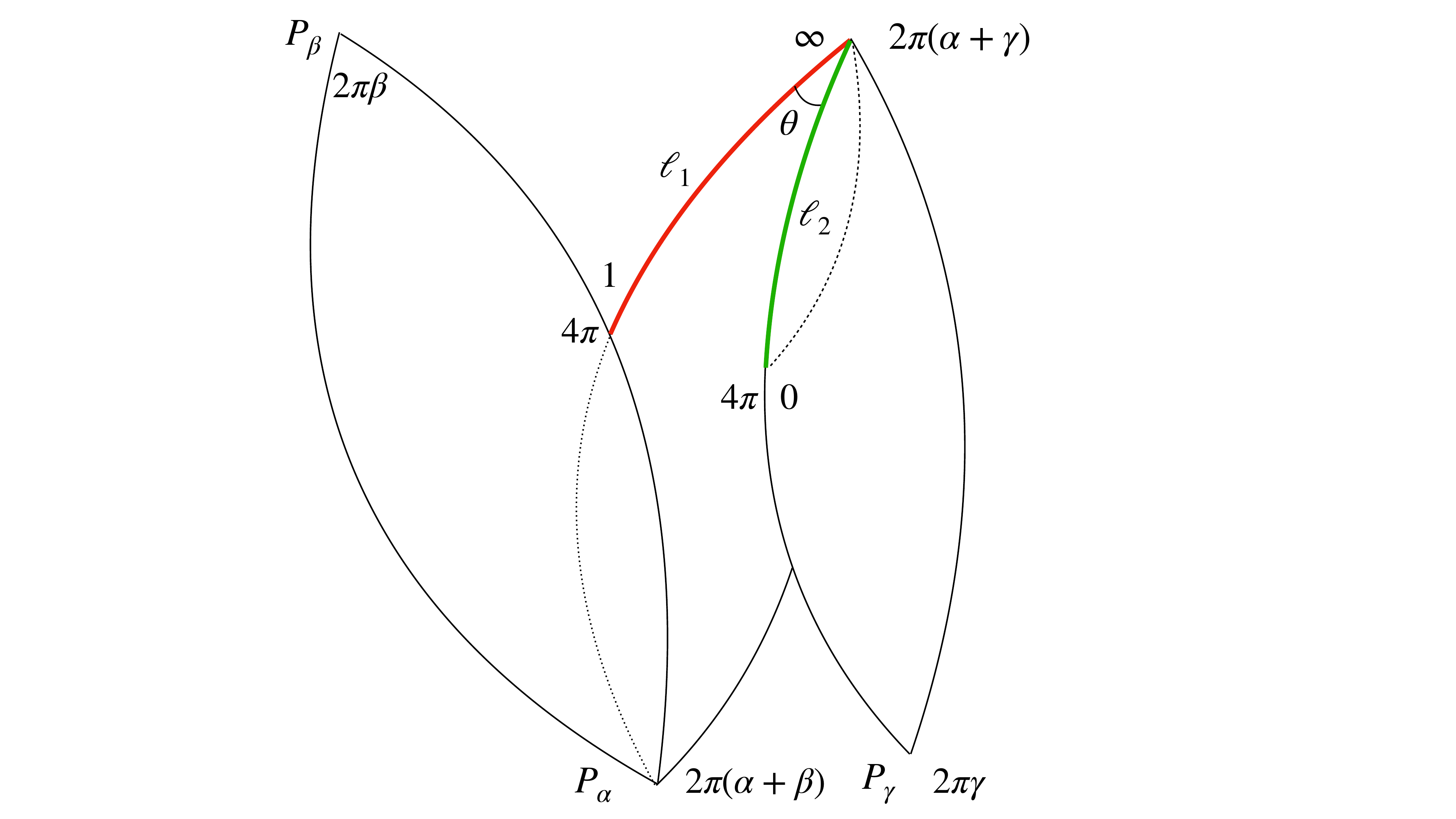} % 
  \caption{A generic configuration of $\bbS^2_{2,2,\beta, \alpha+\beta, \alpha+\gamma, \gamma}$. It can be obtained by gluing three footballs, and is determined by three parameters $\ell_1, \ell_2,\theta$.}
\label{fig:family of threefootbal}
\end{figure}

It turns out that there is another parameter needed to determine the geometry besides $\ell_1, \ell_2$. In  Figure~\ref{fig:threefootball} and \ref{fig:threefootball_v2} we show two distinct configurations which can be seen as the two extreme cases of  Figure~\ref{fig:family of threefootbal}. 
To make up the three real dimensions coming from $c$ and $P_\beta$, we need to add one more parameter. This can be chosen as the angle between two geodesics (see $\theta$ in Figure~\ref{fig:family of threefootbal}.
 The angle can be computed using spherical trigonometry using geodesic lengths from 1 to $\infty$, 0 to $\infty$ and 0 to 1. The first two have been computed above, and the third length $L(0,1)$ is also a function of $c$ and $P_\beta$. It might be possible to find the exact expression by using the explicit metric form~\eqref{eq:ds} and solving for the geodesics, but it is expected to be more complicated than~\eqref{length} so we do not include the computation here, but only note that  $\ell_1, \ell_2, \theta$ locally correspond bijectively to $(c,P_\beta)$. To summarize, we recover the 3-parameter family of geometric decomposition for $\bbS^2_{2,2,\beta, \alpha+\beta, \alpha+\gamma, \gamma}$, which was proved in more general forms in~\cite{Tahar2022} and~\cite{WeiWuXu2022}.

%\begin{equation}
%    L_{total}=2\arctan |F(1)|+2|\arctan |F(0)|-\arctan |F(1)||+ 2\left(\frac{\pi}{2}-\arctan |F(0|\right).
%\end{equation}
%Note that if $\arctan |F(0)| \geq \arctan |F(1)|$, $L_{total}= \pi$; if not, $L=\pi+4\arctan |F(1)|-4\arctan |F(0)|>\pi$. We will show that in case 2 $L_{total}=\pi$ in case 2, while in case 1 the total length is strictly greater than $\pi$. 

\printbibliography 

\end{document}

%% file: prefix.tex
\def\UseBibLatex{1}

\makeatletter
\def\input@path{{styles/}}
\makeatother

\newcommand{\UsePackage}[1]{%
  \IfFileExists{styles/#1.sty}{%
      \usepackage{styles/#1}%
   }{%
      \IfFileExists{../styles/#1.sty}{%
         \usepackage{../styles/#1}%
      }{%
         \usepackage{#1}%
      }%
   }%
}

\usepackage[T1]{fontenc}
\usepackage{lmodern}
\usepackage{textcomp}

\usepackage{amsmath}%
\usepackage{amssymb}%
\usepackage[table]{xcolor}%

\usepackage{todonotes}
\usepackage[in]{fullpage}%

\usepackage[amsmath,thmmarks]{ntheorem}%
\theoremseparator{.}%

\usepackage{titlesec}%
\titlelabel{\thetitle. }%
\usepackage{xcolor}%
\usepackage{mleftright}%
\usepackage{xspace}%
\usepackage{graphicx}
\usepackage{hyperref}%

\usepackage{hyperref}%
\hypersetup{%
      unicode,
      breaklinks,%
      colorlinks=true,%
      urlcolor=[rgb]{0.25,0.0,0.0},%
      linkcolor=[rgb]{0.5,0.0,0.0},%
      citecolor=[rgb]{0,0.2,0.445},%
      filecolor=[rgb]{0,0,0.4},
      anchorcolor=[rgb]={0.0,0.1,0.2}%
}

\theoremseparator{.}%

\theoremstyle{plain}%
\newtheorem{theorem}{Theorem}[section]

\theoremstyle{plain}
\theoremheaderfont{\bfseries}
\theorembodyfont{\upshape}
\newtheorem{remark}{Remark}[section]

\theoremstyle{plain}%

\newcommand{\thmlab}[1]{\label{theo:#1}}

\providecommand{\deflab}[1]{}
\renewcommand{\deflab}[1]{\label{def:#1}}

\newcommand{\remlab}[1]{\label{rem:#1}}

\providecommand{\emphind}[1]{}%
\renewcommand{\emphind}[1]{\emph{#1}\index{#1}}

\definecolor{blue25emph}{rgb}{0, 0, 11}

\providecommand{\emphic}[2]{}
\renewcommand{\emphic}[2]{\textcolor{blue25emph}{%
      \textbf{\emph{#1}}}\index{#2}}

\providecommand{\emphi}[1]{}%
\renewcommand{\emphi}[1]{\emphic{#1}{#1}}

\definecolor{almostblack}{rgb}{0, 0, 0.3}

\providecommand{\emphw}[1]{}%
\renewcommand{\emphw}[1]{{\textcolor{almostblack}{\emph{#1}}}}%

\providecommand{\emphOnly}[1]{}%
\renewcommand{\emphOnly}[1]{\emph{\textcolor{blue25}{\textbf{#1}}}}

\newcommand{\HaoranThanks}[1]{%
   \thanks{Department of Mathematics, Northeastern University, Boston; \href{wu.haoran2@northeastern.edu}{wu.haoran2@northeastern.edu}%
   }%
}

\newcommand{\XuwenThanks}[1]{%
   \thanks{Department of Mathematics, Northeastern University, Boston; \href{x.zhu@northeastern.edu}{x.zhu@northeastern.edu}.
   }
}

\newcommand{\HLink}[2]{\hyperref[#2]{#1~\ref*{#2}}}
\newcommand{\HLinkSuffix}[3]{\hyperref[#2]{#1\ref*{#2}{#3}}}

\providecommand{\deflab}[1]{}
\renewcommand{\deflab}[1]{\label{def:#1}}

\renewcommand{\qed}{\hfill$\square$}
\providecommand{\eqlab}[1]{}%
\renewcommand{\eqlab}[1]{\label{equation:#1}}

\newcommand{\remove}[1]{}%

\usepackage[inline]{enumitem}

\newlist{compactenumA}{enumerate}{5}%
\setlist[compactenumA]{topsep=0pt,itemsep=-1ex,partopsep=1ex,parsep=1ex,%
   label=(\Alph*)}%

\newlist{compactenuma}{enumerate}{5}%
\setlist[compactenuma]{topsep=0pt,itemsep=-1ex,partopsep=1ex,parsep=1ex,%
   label=(\alph*)}%

\newlist{compactenumI}{enumerate}{5}%
\setlist[compactenumI]{topsep=0pt,itemsep=-1ex,partopsep=1ex,parsep=1ex,%
   label=(\Roman*)}%

\newlist{compactenumi}{enumerate}{5}%
\setlist[compactenumi]{topsep=0pt,itemsep=-1ex,partopsep=1ex,parsep=1ex,%
   label=(\roman*)}%

\newlist{compactitem}{itemize}{5}%
\setlist[compactitem]{topsep=0pt,itemsep=-1ex,partopsep=1ex,parsep=1ex,%
   label=\ensuremath{\bullet}}%

\providecommand{\BibLatexMode}[1]{}
\providecommand{\BibTexMode}[1]{}

\ifx\UseBibLatex\undefined%
  \renewcommand{\BibLatexMode}[1]{}
  \renewcommand{\BibTexMode}[1]{#1}
\else
  \renewcommand{\BibLatexMode}[1]{#1}
  \renewcommand{\BibTexMode}[1]{}
\fi

\BibLatexMode{%
   \usepackage[bibencoding=utf8,style=numeric,citestyle=numeric,backend=biber,sorting=nyt,]{biblatex}%
   \UsePackage{my_biblatex}%
}

\numberwithin{figure}{section}%
\numberwithin{table}{section}%
\numberwithin{equation}{section}%